\documentclass[a4paper,12pt]{amsart}

\usepackage{amsfonts}

\def\CP{{\mathbb{C}\mathrm{P}^1}}
\def\aut{{\mathrm{aut}}}
\def\oM{{\overline{\mathcal{M}}}}
\def\M{{{\mathcal{M}}}}
\def\d{{\partial}}

\newtheorem{theorem}{Theorem}
\newtheorem{lemma}{Lemma}

\title{Some relations for one-part double Hurwitz numbers}

\thanks{This work is partially supported by RFBR
grant no. 02-01-22004.}

\author{S.~V.~Shadrin}

\date{}

\begin{document}

\maketitle

\textbf{1.}
In this note we state some new relations for double Hurwitz
numbers. These relations come from the intersection theory of
the moduli spaces of curves and generalize the results announced
in~\cite{s1} and proved in~\cite{s2}.

The Hurwitz numbers studied here have recently been considered
in details in~\cite{gjv}. In particular, we refer
to~\cite{gjv} as to a good survey in this area with the most
complete collection of references.

\textbf{2.} Recall the definition of Hurwitz numbers. Fix
integers $g\geq 0$, $n\geq 1$, and two unordered partitions of
the number $n$, $n=a_1+\dots+a_p=b_1+\dots+b_q$.
Consider $\CP$ with marked $0$, $\infty$, and
$2g+p+q-2$ points more, $z_1,\dots,z_{2g+p+q-2}$.

There is a finite number of functions $f\colon C\to\CP$ of
degree $n$ defined on curves of genus $g$ such that $(f) =
-\sum_{i=1}^p a_i x_i + \sum_{j=1}^q b_j y_j$ (here $x_i,y_j$
are pairwise distinct points of $C$), and
$z_1,\dots,z_{2g+p+q-2}$ are simple critical values of $f$.

The number of such functions is called \emph{Hurwitz number} and
denoted by $H_g(a_1,\dots,a_p|b_1,\dots,b_q)$. Here we count
each function $f\colon C\to\CP$ weighted by $1/|\aut(f)|$.
For instance, $H_g(2|2)=1/2$ for any $g$.

We will consider only the numbers $H_g(n|b_1,\dots,b_q)$.
For convenience we introduce the following notation:
\begin{equation}
\widehat H_{g,n}(b_1,\dots,b_q) := \frac
{|\aut(b_1,\dots,b_q)|} {n^{2g+q-2}(2g+q-1)!}
H_g(n|b_1,\dots,b_q).
\end{equation}

\textbf{3.} Let us state the main theorem:
\begin{theorem} \label{rel}
For any $g\geq 0$, $n\geq 1$, $b_1+\dots+b_q=n$, we have:
\begin{multline}
\binom{g}{0} \widehat H_{g,n}(b_1,\dots,b_q) -
\binom{g}{1} \widehat H_{g,n+1}(b_1,\dots,b_q,1) + \\
+ \dots
+ (-1)^g \binom{g}{g} \widehat
H_{g,n+g}(b_1,\dots,b_q,\underbrace{1,\dots,1}_g) =
\frac {(-1)^g} {24^g}.
\end{multline}
\end{theorem}

Let us check this theorem in a special case. Let $g=2$.
Then there is a fomula from~\cite{gjv}:
\begin{equation}
\widehat H_{g,n}(b_1,\dots,b_q)=\frac{1}{24^2}\left(
\frac{1}{2}A^2-\frac{1}{5}B \right),
\end{equation}
where $A=-1+b_1^2+\dots+b_q^2$, $B=-1+b_1^4+\dots+b_q^4$.

Therefore, in this special case our theorem is equivalent to the
following equation:
\begin{multline}
\left( \frac{1}{2}A^2-\frac{1}{5}B \right)
-2\left( \frac{1}{2}(A+1)^2-\frac{1}{5}(B+1) \right) + \\
+\left( \frac{1}{2}(A+2)^2-\frac{1}{5}(B+2) \right)=1.
\end{multline}
One can check this by direct calculation.

\textbf{4.} Let us explain, how one can obtain such relation for
Hurwitz numbers.

Consider the moduli space $\oM_{g,q+1}$ of curves of genus $g$
with $q+1$ marked points. Let $L_1$ be the line bundle over
$\oM_{g,q+1}$ with the fiber $T^*_{x_1}C$ at a point
$(C,x_1,\dots,x_{q+1})\in\oM_{g,q+1}$. By $\psi_1$ denote
$c_1(L_1)$ .

Consider the subvariety $V^\circ_g(b_1,\dots,b_q)\subset
\M_{g,q+1}$ consisting of smooth curves $(C,x_1,\dots,x_{q+1})$
such that $-(\sum_{i=1}^qb_i)x_1+b_1x_2+\dots+b_qx_{q+1}$ is a
divisor of a meromorphic function. By $V_g(b_1,\dots,b_q)$
denote the closure of $V^\circ_g(b_1,\dots,b_q)$ in
$\oM_{g,q+1}$.

Theorem~\ref{rel} is the direct corollary of two
relations for intersection numbers on the moduli spaces of
curves looking like follows:

\begin{lemma} For any $b_1,\dots,b_q$,
\begin{multline}
\frac{(-1)^g}{24^g}=(-1)^g g!\int_{\oM_{g,3}}\psi_1^{3g}=
\binom{g}{0}\int_{V_g(b_1,\dots,b_q)}\psi_1^{2g+q-2}
-
\\
\binom{g}{1}\int_{V_g(b_1,\dots,b_q,1)}\psi_1^{2g+q-1}
+ \dots +(-1)^g
\binom{g}{g}\int_{V_g(b_1,\dots,b_q,1,\dots,1)}\psi_1^{3g+q-2}.
\end{multline}
\end{lemma}

\begin{lemma} For any $b_1,\dots,b_q$,
\begin{equation}
\widehat H_{g,n}(b_1,\dots,b_q) =
 \int_{V_g(b_1,\dots,b_q)}\psi_1^{2g+q-2}
\end{equation}
\end{lemma}

Lemma 1 is a generalization of the similar statement
in~\cite{s1,s2}. Lemma 2 is just one of the theorems
proved in~\cite{s2}.

\textbf{5.} There is a `cut-and-join' type equation for
the generating function of Hurwitz numbers considered here.
(for the similar equations for some other Hurwitz numbers,
see~\cite{gjvn, z}).

Consider the function $F(\theta,x_1,x_2,\dots)$ defined like
follows:
\begin{multline}
F(\theta,x_1,x_2,\dots)=\\
=\sum_{g\geq 0}\sum_{(b_1,\dots,b_q)}
H_g(n|b_1,\dots,b_q)
\frac{\theta^{(2g+q-1)}}
{(2g+q-1)!}x_{b_1}\dots x_{b_q}.
\end{multline}

\begin{theorem}
The function $F=F(\theta,x_1,x_2,\dots)$ satisfies
the following:
\begin{equation}
\frac{\d F}{\d \theta}=\frac{1}{2}
\sum_{i,j\geq 1}\left(
ijx_{i+j}\frac{\d^2 F}{\d x_i \d x_j}
+(i+j)x_ix_j\frac{\d F}{\d x_{i+j}}
\right).
\end{equation}
\end{theorem}

Generally speaking, this theorem is more or less obvious
from the point of view of geometry of ramified coverings.
Surely, this theorem has the similar purely combinatorial proof
as its analog in~\cite{gjvn}. But this theorem is also a direct
corollary of relations for intersection numbers obtained
in~\cite{s2}.

Really, one can rewrite this theorem as a relation for Hurwitz
numbers like follows:

\begin{multline}
H'_g(n|b_1,\dots,b_q)=\\
\frac{1}{2}
\sum_{i=1}^q\sum_{b_i^{(1)}+b_i^{(2)}=b_i}b_i^{(1)}b_i^{(2)}
H'_{g-1}(n|b_1,\dots,\widehat{b_i},\dots,b_q,b_i^{(1)},
b_i^{(2)}) + \\
\frac{1}{2}
\sum_{i\not= j}
(b_i+b_j)
H'_{g}(n|b_1,\dots,\widehat{b_i},\dots,\widehat{b_j},
\dots,b_q, b_i+b_j),
\end{multline}
where by $H'_g(n|b_1,\dots,b_q)$ we denote
$$
H_g(n|b_1,\dots,b_q)\cdot |\aut(b_1,\dots,b_q)|=
\widehat H_{g,n}(b_1,\dots,b_q) \cdot
n^{2g+q-2}(2g+q-1)!
$$

Using Lemma 2, we can rewrite Eqn.~(9) as a relation
for the integrals
$\int_{V_g(b_1,\dots,b_q)}\psi_1^{2g+q-2}$.
Thus we obtain exactly a special case
of Theorem 12.2 from~\cite{s2}.

\bigskip

\noindent
Independent University of Moscow\\
Stockholm University\\
e-mail: shadrin@mccme.ru


\begin{thebibliography}{0}

\bibitem{gjvn} I.~P.~Goulden, D.~M.~Jackson, A.~Vainshtein, Ann.
Comb. 4 (2000),27-46

\bibitem{gjv} I.~P.~Goulden, D.~M.~Jackson, R.~Vakil, arXiv:
math.AG/0309440

\bibitem{s1} S.~V.~Shadrin, Russ.~Math.~Surv, 58:1, 195-196
(2003)

\bibitem{s2} S.~V.~Shadrin, Int.~Math.~Res.~Not. 2003, no. 38,
2051-2094

\bibitem{z} D. Zvonkine, arXiv: math.AG/0304251

\end{thebibliography}
\end{document}